\documentclass[12pt]{article}
\usepackage[reqno]{amsmath}
\usepackage{amssymb, theorem, enumerate}
\usepackage{graphicx,psfrag}
\usepackage{pgf,tikz,pgfplots}
\usetikzlibrary{automata, positioning, arrows}
\usepackage[abbrev,msc-links]{amsrefs}
\usetikzlibrary{decorations.markings}
\usepackage{array}
\usepackage{float}
\usepackage{wrapfig}
\usepackage{fullpage}
\usepackage{url}
\usepackage{mathtools}
\pgfplotsset{compat=1.18}

\theoremstyle{change}
{\theorembodyfont{\slshape}
\newtheorem{theorem}{Theorem.}[section]
\newtheorem{lemma}[theorem]{Lemma.}
\newtheorem{remark}[theorem]{Remark.}
\newtheorem{corollary}[theorem]{Corollary.}}

\newtheorem{problem}[theorem]{Problem.}
\theorembodyfont{\rmfamily}
\newtheorem{eg}[theorem]{Example.}

\newtheorem{Local switching}[theorem]{Local switching.}

\newcommand\cref[1]{Corollary~\ref{cor:#1}}

\def\proof{\noindent{{\sl Proof. }}}
\def\sqr#1#2{{\vbox{\hrule height.#2pt
    \hbox{\vrule width.#2pt height#1pt \kern#1pt
        \vrule width.#2pt}\hrule height.#2pt}}}
\def\eqed{\sqr53}
\def\qed{%
    \ifmmode\eqno\eqed
    \else\nobreak\ \hfill\eqed\medbreak\fi}

\newcommand\al{\alpha}

\newcommand\ga{\gamma}

\newcommand\fld{{\mathbb F}}

\newcommand\re{{\mathbb R}}
\newcommand\rats{{\mathbb Q}}

\newcommand\comp[1]{{\mkern2mu\overline{\mkern-2mu#1}}}

\newcommand\pmat[1]{\begin{pmatrix} #1 \end{pmatrix}}

\DeclareMathOperator{\diag}{diag}

\newcommand{\ptmu}[1]{\psi(#1,t,\mu)}

\title{Degree-Similar Graphs}

\author{Chris Godsil$^1$\footnote{C. Godsil gratefully acknowledges the support of the
Natural Sciences and Engineering Council of Canada (NSERC), Grant No. RGPIN-9439.},
Wanting Sun$^{1,2}$\footnote{Most of this research was done when W. Sun was a visiting Ph. D. student at University of Waterloo. W. Sun is  supported by the China Postdoctoral Science Foundation (2023M742092).\newline
\hspace*{5mm}{\it Email addresses}: cgodsil@uwaterloo.ca (C. Godsil),\ wtsun2018@sina.com (W. Sun).}\\
 	$^1$Department of Combinatorics \& Optimization, University of\\
	Waterloo, Waterloo, Ontario, Canada\\
$^2$Data Science Institute, Shandong, Jinan, 250100, PR China}

\begin{document}
\maketitle

\begin{abstract}
	The degree matrix of a graph is the diagonal matrix with diagonal entries
	equal to the degrees of the vertices of $X$. If $X_1$ and $X_2$ are graphs
	with respective adjacency matrices $A_1$
	and $A_2$ and degree matrices $D_1$ and $D_2$, we say that $X_1$ and $X_2$
	are \textsl{degree similar} if there is an invertible real matrix $M$ such
	that $M^{-1}A_1M=A_2$ and $M^{-1}D_1M=D_2$. If graphs $X_1$ and $X_2$ are
	degree similar, then their adjacency matrices, Laplacian matrices, unsigned
	Laplacian matrices and normalized Laplacian matrices are similar. We first show
	that the converse is not true. Then, we provide a number of constructions of
	degree-similar graphs. Finally, we show that the matrices $A_1-\mu D_1$ and
	$A_2-\mu D_2$ are similar over the field of rational functions $\rats(\mu)$
	if and only if the Smith normal forms of the matrices $tI-(A_1-\mu D_1)$
	and $tI-(A_2-\mu D_2)$ are equal.
\end{abstract}

\section{Introduction}

Let $X$ be a graph with vertex set $V(X)$ and edge set $E(X)$. We use $A=A(X)$ to denote the \textsl{adjacency matrix} of $X$ and $D=D(X)$ to denote
the \textsl{degree matrix}, the diagonal matrix
with $D_{i,i}$ equals to the valency of the vertex $i$ in $X$.
If $X_1$ and $X_2$ are graphs with respective adjacency matrices $A_1$
and $A_2$ and degree matrices $D_1$ and $D_2$, we say that $X_1$ and $X_2$
are \textsl{degree similar} if there is an invertible real matrix $M$ such that
\begin{equation}\label{eq:3}
	M^{-1}A_1M= A_2,\quad M^{-1}D_1M = D_2.
\end{equation}
Clearly, if $X_1$ and $X_2$ are degree similar, then their adjacency matrices,
Laplacians $D-A$, unsigned Laplacians $D+A$, and their normalized Laplacians
$D^{-1/2}AD^{-1/2}$ are similar. (When using the normalized Laplacian,
we assume the underlying
graph has no isolated vertices.) Thus we have a hierarchy of conditions
on a pair of graphs:
\begin{enumerate}[(a)]
	\item
	They are degree similar.
	\item
	Their adjacency matrices and the three Laplacians are similar.
	\item
	Their adjacency matrices are similar.
\end{enumerate}
We note that these three conditions are equivalent for regular graphs.

We discuss some related earlier work.
Butler et al. \cite{butler2022complements} constructed graphs that are
cospectral with respect to adjacency, Laplacian,
unsigned Laplacian and normalized Laplacian matrices.
In \cite{wang2011graphs}, by using local switching, Wang et al. gave a construction of pairs of degree-similar graphs.
Guo et al. \cite{guo2013laplacian}  derived six reduction procedures on the Laplacian, unsigned Laplacian
and normalized Laplacian characteristic polynomials of a graph which can be
used to construct larger Laplacian, unsigned Laplacian and normalized Laplacian
cospectral graphs, respectively.

Tutte \cite{tutte1979all} defined the \textsl{idiosyncratic polynomial} of
a graph to be 
\[ p(x,\al):=\det(A + \al(J-I-A) - xI),\] 
where $I$ is the identity matrix and $J$ is the all 1s matrix of appropriate size.
Van Dam and Haemers \cite{vanDam} worked with
the \textsl{generalized adjacency matrix}
$sI+tA+\mu J$, the characteristic polynomial of this  matrix is a form of the idiosyncratic polynomial.

Wang et al.~\cite{wang2011graphs} defined the generalized
characteristic polynomial $\psi(X,t,\mu)$ of a graph $X$ as follows:
\[
	\psi(X,t,\mu) := \det(tI - (A-\mu D).
\]
Note that if $\ptmu{X_1}=\ptmu{X_2}$, the adjacency matrices of $X_1$
and $X_2$ are similar, along with the three Laplacians.
They observed that if $X_1$ and $X_2$ are degree similar, then
$\ptmu{X_1}=\ptmu{X_2}$, and asked if the converse was true.
Our results in Section~3 show that it is not. Hence if
the adjacency and the three Laplacian matrices
of two graphs are similar, it does not follow that the graphs are
degree similar.

In Section 4, we prove that if $X$ and $Y$ are connected graphs and are degree similar,
then their complements $\comp{X}$ and $\comp{Y}$ are degree similar.
In Sections 5-8, we provide a number of constructions of pairs of (non-isomorphic)
degree-similar graphs, for example, graph products, adding or deleting vertices and so on.
In Section 9, we study the relation between similarity and Smith normal forms of matrices. In the last section, we provide some further discussions.

\section{Ihara zeta function}

In this section, we note one further consequence of degree similarity.

A walk in a graph $X$ is \textsl{reduced} if it does not contain
any subsequence of the form $uvu$; such walks may also be called
\textsl{non-backtracking}. If $|V(X)|=n$, then $p_r(A)$
denotes the $n\times n$ matrix where $(p_r(A))_{u,v}$ is the number
of reduced walks in $X$ from $u$ to $v$. So
\[
	p_0(A)=I,\ p_1(A) = A,\ p_2(A) = A^2-D.
\]
When $r\ge3$, we have the recurrence
\[
	Ap_r(A) = p_{r+1}(A) + (D-I)p_{r-1}(A),
\]
from which it follows that $p_r(A)$ is a polynomial in $A$ and $D$.
These observations are due to Biggs. For details, and
for the following theorem, see Chan and Godsil \cite{chan1997symmetry}.

\begin{theorem}
	For any connected graph on at least two vertices,
	\[
		\sum_{r\ge0} t^r p_r(A) = (1-t^2)(I-tA+t^2(D-I))^{-1}.
	\]
\end{theorem}

The determinant of the generating function on the left in this identity is the
\textsl{Ihara zeta function} of the graph, and therefore
if $X_1$ and $X_2$ are degree similar and have no isolated
vertices, their Ihara zeta functions are equal. (This fact was noted by
Wang et al.~in \cite{wang2011graphs}, with a sketch of a proof. For more on
Ihara zeta functions, see \cite{terras2010zeta}.)

\section{Trees}

Firstly, we describe some notation. For a graph $X$ and a vertex $u\in V(X)$, we use $d_X(u)$ to denote the degree of $u$ in $X$. For a vertex subset
$U\subseteq V(X)$, denote the induced subgraph of $X$ on $U$ by $X[U]$, and the induced subgraph of $X$ on $V(X)\setminus U$ by $X\backslash U$. 
For an $n\times n$ matrix $M$ and a set $U\subset \{1,\ldots,n\}$, we use $M(U)$
to denote a matrix obtained from $M$ by deleting the rows in $U$
and the columns in $U$.
When $U=\{u\}$, we use $X\backslash u$ and $M(u)$ instead.
If ${\bf x}$ and ${\bf y}$ are
two column vectors, we use $[M|{\bf x},{\bf y}]$ to denote the bordered matrix
\[
	\pmat{0&{\bf x}^T\\{\bf y}&M}.
\]

Assume that $S$ is a graph and $T$ is a rooted tree.
The \textsl{coalescence} $S\bullet T$
is the graph formed by identifying the root of $T$ and a vertex of $S$.

\begin{figure}[ht!]
  \centering
  \begin{tikzpicture}[scale = 1.5]
  \tikzstyle{vertex}=[circle,fill=black,minimum size=0.5em,inner sep=0pt]
  \node[vertex] (G_1) at (0,0)[label=right:$v$]{};%
  \node[vertex] (G_2) at (-0.3,0.5){};
  \node[vertex] (G_3) at (-0.6,1){};
  \node[vertex] (G_4) at (-0.9,1.5){};
  \node[vertex] (G_5) at (-0.05,1){};
  \node[vertex] (G_6) at (0.2,1.5){};
  \node[vertex] (G_7) at (0.3,0.5){};
  \node[vertex] (G_8) at (0.6,1){};
  \node[vertex] (G_9) at (0.4,1.5){};
  \node[vertex] (G_10) at (0.2,2){};
  \node[vertex] (G_11) at (0,2.5){};
  \node[vertex] (G_12) at (0.9,1.5){};
  \node[vertex] (G_13) at (1.2,2){};
  \node[vertex] (G_14) at (1.5,2.5){};
  \node[vertex] (G_15) at (1.2,3){};
  \node[vertex] (G_16) at (1.8,3){};
  \draw[thick] (G_1)--(G_2)--(G_3)--(G_4);
  \draw[thick] (G_2)--(G_5)--(G_6);
  \draw[thick] (G_1)--(G_7) -- (G_8)--(G_9)--(G_10)--(G_11);
  \draw[thick] (G_8)--(G_12)--(G_13)--(G_14)--(G_15);
  \draw[thick] (G_14)--(G_16);
  \draw (0,-0.4)node{$T_1$};
  \end{tikzpicture}
\hspace{2em}
  \begin{tikzpicture}[scale = 1.5]
  \tikzstyle{vertex}=[circle,fill=black,minimum size=0.5em,inner sep=0pt]
  \node[vertex] (G_1) at (0,0)[label=right:$w$]{};%
  \node[vertex] (G_2) at (-0.3,0.5){};
  \node[vertex] (G_3) at (-0.6,1){};
  \node[vertex] (G_4) at (-0.9,1.5){};
  \node[vertex] (G_5) at (-0.4,1.5){};
  \node[vertex] (G_6) at (0.3,0.5){};
  \node[vertex] (G_7) at (0.05,1){};
  \node[vertex] (G_8) at (-0.2,1.5){};
  \node[vertex] (G_9) at (-0.45,2){};
  \node[vertex] (G_10) at (0.6,1){};
  \node[vertex] (G_11) at (0.9,1.5){};
  \node[vertex] (G_12) at (1.2,2){};
  \node[vertex] (G_13) at (0.9,2.5){};
  \node[vertex] (G_14) at (0.6,3){};
  \node[vertex] (G_15) at (1.5,2.5){};
  \node[vertex] (G_16) at (1.8,3){};
  \draw[thick] (G_1)--(G_2)--(G_3)--(G_4);
  \draw[thick] (G_3)--(G_5);
  \draw[thick] (G_1)--(G_6)--(G_7) -- (G_8)--(G_9);
  \draw[thick] (G_6)--(G_10)--(G_11)--(G_12)--(G_13)--(G_14);
  \draw[thick] (G_12)--(G_15)--(G_16);
  \draw (0,-0.4)node{$T_2$};
  \end{tikzpicture}
  \caption{$T_1$ and $T_2$.}\label{fig001}
\end{figure}
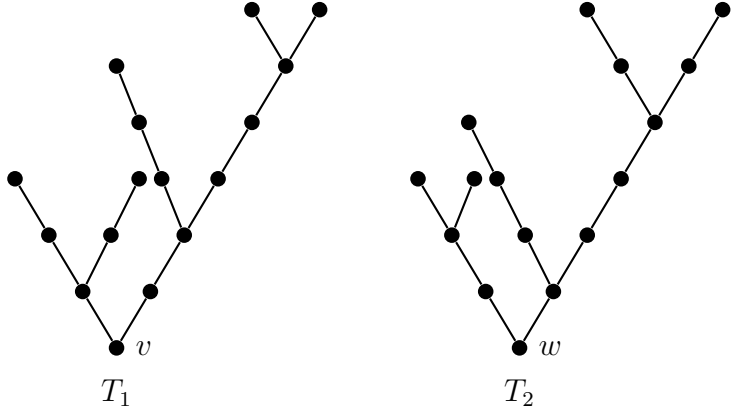

Let $T_1$ and $T_2$ be the rooted trees shown in Figure \ref{fig001}, whose roots are $v$ and $w$ respectively. Clearly, $T_1$
and $T_2$ are two isomorphic trees with different roots.
In 1977, McKay \cite{mckay1977spectral}
showed that for any tree $S$ with at least two vertices, $S\bullet T_1$
and $S\bullet T_2$ are not isomorphic, but they are cospectral with respect to
several graph matrices (in particular, adjacency matrix, Laplacian matrix and
unsigned Laplacian matrix). Osborne \cite{osborne2013} showed that their normalized Laplacian matrices are also similar.

In fact, we can prove a more general result: for any graph $S$ with at least two
vertices, $\psi(S\bullet T_1,t,\mu)=\psi(S\bullet T_2,t,\mu)$,
i.e., $S\bullet T_1$ and $S\bullet T_2$ are cospectral with respect to
the matrix $A-\mu D$. For convenience, put $A_\mu(X):=A(X)-\mu D(X)$.

\begin{lemma}\label{lem7}
    Let $T_1$ and $T_2$ be the two trees depicted in Figure \ref{fig001}, and let $S_i=S\bullet T_i$ for $i=1,2$, where $S$ is a non-trivial
    graph. Then $\psi(S_1,t,\mu)=\psi(S_2,t,\mu)$.
\end{lemma}

\proof
  Without loss of generality, assume that $S_i$ is obtained by identifying the root of $T_i$ and a vertex $r$ of $S$, here $i=1,2$. It is routine to check that
    \begin{equation*}
      A_\mu(S_1)=\left(
      \begin{matrix}
        -\mu(d_S(r)+2) & {\bf x}_{T_1}^T & {\bf y}_S^T \\
        {\bf x}_{T_1} & A_\mu (T_1)(v) & {\bf 0} \\
        {\bf y}_S & {\bf 0} & A_\mu (S)(r) \\
      \end{matrix}
    \right)
    \end{equation*}
  for some column vectors ${\bf x}_{T_1}$ and ${\bf y}_S$.  For convenience,
  denote $\phi(M):=\det(tI-M)$. Based on \cite{mckay1977spectral}*{Lemma 2.2(i)},
  one has
    \begin{align*}
    \psi(S_1,t,\mu)=&\phi(A_\mu (T_1)(v))\phi([A_\mu(S)(r)|{\bf y}_S,{\bf y}_S])+\phi(A_\mu(S)(r))\phi([A_\mu (T_1)(v)|{\bf x}_{T_1},{\bf x}_{T_1}])\\
    	&-(t-\mu(d_S(r)+2))\phi(A_\mu (T_1)(v))\phi(A_\mu(S)(r)).
    \end{align*}
    Similarly, one may write $\psi(S_2,t,\mu)$ as follows:
    \begin{align*}
    \psi(S_2,t,\mu)=&\phi(A_\mu (T_2)(w))\phi([A_\mu(S)(r)|{\bf y}_S,{\bf y}_S])+\phi(A_\mu(S)(r))\phi([A_\mu (T_2)(w)|{\bf x}_{T_2},{\bf x}_{T_2}])\\
    	&-(t-\mu(d_S(r)+2))\phi(A_\mu (T_2)(w))\phi(A_\mu(S)(r)).
    \end{align*}
    By a direct calculation, we obtain
    \[
    \phi(A_\mu (T_1)(v))=\phi(A_\mu (T_2)(w))
    \]
     and
     \[
     \phi([A_\mu (T_1)(v)|{\bf x}_{T_1},{\bf x}_{T_1}])
     	=\phi([A_\mu (T_2)(w)|{\bf x}_{T_2},{\bf x}_{T_2}]).
      \]
     It follows that $\psi(S_1,t,\mu)=\psi(S_2,t,\mu)$.\qed

Obviously, if $X_1$ and $X_2$ are degree similar,
then $\psi(X_1,t,\mu)=\psi(X_2,t,\mu)$.
Wang et al. \cite{wang2011graphs} proposed a problem: Is the converse true? Next, we give
some examples to show that it is not.

The following result is a reformulation of \cite{mckay1977spectral}*{Theorem 5.3}.

\begin{theorem}\label{lem6}
  Two trees are degree similar if and only if they are isomorphic.
\end{theorem}

Combining Lemma \ref{lem7} with Theorem \ref{lem6}, we have the following corollary.

\begin{corollary}
\label{cor:notdgs}
For any tree $S$ with at least two vertices, we have $\psi(S\bullet T_1,t,\mu)=\psi(S\bullet T_2,t,\mu)$,
but $S\bullet T_1$ and $S\bullet T_2$ are not degree similar.
\end{corollary}

\section{Subgraphs and complements}

In this section, we study the subgraphs and complements of degree-similar graphs.
Recalling the definition of degree-similar graphs, we first present a
basic property of the invertible real matrix $M$ in \eqref{eq:3}.

\begin{lemma}\label{lem1}
  Let $X_1$ and $X_2$ be graphs with degree matrices $D_1$ and $D_2$ respectively.
If there is an invertible real matrix $M$ such that
\[
	M^{-1}D_1M = D_2.
\]
Then $M$ is block diagonal.
\end{lemma}

\proof
  Assume that $d_1,\ldots,d_t$ are all different vertex degrees of $X_1$.
  Partition the vertex set of $X_1$ as follows:
  $V(X_1)=V_1\cup V_2\cup \cdots \cup V_t$, where $V_i=\{w:d_{X_1}(w)=d_i\}$ for $i\in \{1,\ldots,t\}$.
  Since $D_1$ is a diagonal matrix, after reordering the vertices of $X_1$,
  we can write $D_1$ as follows:
  \begin{equation*}
    D_1=\left(
          \begin{array}{cccc}
            d_1I_{|V_1|} &  &  &  \\
             & d_2I_{|V_2|} &  &  \\
             &  & \ddots &  \\
             &  &  & d_tI_{|V_t|} \\
          \end{array}
        \right).
  \end{equation*}
  Notice that $D_1$ and $D_2$ are diagonal matrices.
  Together with $M^{-1}D_1M = D_2$, there exists a
  permutation matrix $P$ such that $P^TM^{-1}D_1MP = P^TD_2P=D_1$.
  Let $Q=MP$. Then $Q^{-1}D_1Q=D_1$, i.e., $D_1Q=QD_1$.
  Therefore, $Q$ is a block diagonal matrix with respect to the
  partition $V_1\cup V_2\cup \cdots \cup V_t$, which implies that $M$
  is block diagonal.\qed

The following result is an immediate consequence of Lemma~\ref{lem1},
which gives some cospectral graphs with respect to the adjacency matrix.

\begin{lemma}\label{lem2}
	Let $X_1$ and $X_2$ be two degree-similar graphs, and let $d$ be the
	degree of some vertex in $X_1$. Assume that $V_i=\{w:d_{X_i}(w)=d\}$
	for $i\in \{1,2\}$. Then the induced subgraphs $X_1[V_1]$ and $X_2[V_2]$
	are adjacency cospectral.
\end{lemma}

\begin{remark}
  In fact, Lemma \ref{lem2} is more useful in determining two graphs that are
  not degree similar. For example, let $X$ be a strongly regular graph with
  parameters SRG$(25,12, 5, 6)$. In fact, there are exactly $15$ non-isomorphic
  strongly regular graphs with such parameters. Here, we assume the adjacency
  matrix of $X$ is the first one described in Spence's website:
\url{http://www.maths.gla.ac.uk/~es/srgraphs.php}.

Assume that the first two rows of $A(X)$ are indexed by $u$ and $v$
respectively. One may check
  \begin{align*}
    &\det(xI_{12}-A(X\backslash N_X[u]))-\det(xI_{12}-A(X\backslash N_X[v]))\\
    =&-2x^9+2x^8+64x^7+39x^6-372x^5-135x^4+648x^3-324x^2,
  \end{align*}
  here $N_X[u]$ denotes the closed neighborhood of $u$ in $X$. This implies
  that $X\backslash N_X[u]$ and $X\backslash N_X[v]$ are not adjacency cospectral.
  Together with Lemma \ref{lem2},  we know $X\backslash u$
  and $X\backslash v$ are not degree similar.\qed
\end{remark}

Next, we show that degree similar is preserved under taking the complement
of the underlying graphs. Let $\overline{X}$ denote the complement of a
graph $X$.

\begin{lemma}\label{lem4}
  If $X$ is connected, $X$ and $Y$ are degree-similar, then their
  complements are degree similar.
\end{lemma}

\proof
Assume that $X$ and $Y$ have $n$ vertices. Since $X$ and $Y$ are degree similar,
there exists an invertible real matrix $M$ such that
\[
	M^{-1}A(X)M= A(Y),\quad M^{-1}D(X)M = D(Y).
\]
Since the Laplacians of $X$ and $Y$ are cospectral, $X_2$ is
connected. Then there is a polynomial $p$, determined by the
spectrum of $D(X)-A(X)$, such that $p(D(X)-A(X))=J_n$. Therefore $p(D(Y)-A(Y))=J_n$. Consequently,
\begin{equation}\label{eq:1}
	M^{-1}J_nM = M^{-1}p(D(X)-A(X))M = p(D(Y)-A(Y)) = J_n,
\end{equation}
from which it follows that $J_n-I_n-A(X)$ and $J_n-I_n-A(Y)$ are cospectral.

Notice that $A(\overline{X})=J_n-I_n-A(X)$ and $D(\overline{X})=(n-1)I_n-D(X)$. Then,
\begin{align*}
  M^{-1}A(\overline{X})M&=M^{-1}(J_n-I_n-A(X))M=J_n-I_n-A(Y)=A(\overline{Y}),\\
   M^{-1}D(\overline{X})M&=M^{-1}((n-1)I_n-D(X))M=(n-1)I_n-D(Y)=D(\overline{Y}).
\end{align*}
It follows that $\overline{X}$ and $\overline{Y}$ are degree similar. \qed

\section{Local switching}

There is a powerful and productive method called
\textsl{local switching} \cite{GM},
which can produce numerous pairs of cospectral graphs.
Wang et al.~\cite{wang2011graphs}
constructed a family of degree-similar graphs by using local switching.
In this section, we generalize their result, and use local switching to
construct a large family of degree-similar graphs.
Firstly, we describe local switching.

\begin{Local switching}
  Let $X$ be a graph and let $\pi:=C_1\cup C_2\cup \cdots \cup C_k\cup C$
    be a partition of $V(X)$.
  Suppose that, whenever $1\leqslant i,j\leqslant k$ and $v\in C$, we have
  \begin{enumerate}[(a)]
    \item any two vertices in $C_i$ have the same number of
    neighbors in $C_j$, and
    \item $v$ has either $0$, $\frac{|C_i|}{2}$ or $|C_i|$ neighbors in $C_i$.
  \end{enumerate}
  The graph $X^{\pi}$ formed by local switching in $X$ with respect to $\pi$ is
  obtained from $X$ as follows. For each $v\in C$ and $1\leqslant i\leqslant k$
  such that $v$ has $\frac{|C_i|}{2}$ neighbors in $C_i$, delete
  those $\frac{|C_i|}{2}$ edges and join $v$ instead to the other $\frac{|C_i|}{2}$ vertices in $C_i$.
\end{Local switching}

Godsil and McKay \cite{GM} showed that if $X^{\pi}$ is the graph formed by
local switching in $X$ with respect to a partition $\pi$, then $X$ and $X^{\pi}$
are cospectral, with cospectral complements. Now, we use local switching to
construct a large family of degree-similar graphs.

\begin{lemma}\label{lem9}
  Let $X$ be a graph and let $\pi:=C_1\cup C_2\cup \cdots \cup C_k\cup C$ be a partition
  of $V(X)$. Put $c:=|C|$ and $c_i:=|C_i|$ for $1\leqslant i\leqslant k$.
  Suppose that, whenever $1\leqslant i,j\leqslant k$, we have
  \begin{enumerate}[(a)]
    \item any two vertices in $C_i$ have $d_{ij}$ neighbors in $C_j$,
    \item for any $u\in C$, $u$ has $\frac{c_i}{2}$ neighbors in $C_i$,  and
    \item for any $v\in C_i$, $v$ has $\frac{c}{2}$ neighbors in $C$.
  \end{enumerate}
  If $X^{\pi}$ is formed by local switching in $X$ with respect to $\pi$,
  then $X$ and $X^{\pi}$ are degree similar.
\end{lemma}

\proof
Assume that the vertices of $X$ are labelled in an order consistent with $\pi$.
For $i\in \{1,\ldots,k\}$, let $Q_i=\frac{2}{c_i}J_{c_i}-I_{c_i}$.
Clearly, $Q_i$ is an orthogonal matrix and $Q_i=Q_i^T=Q_i^{-1}$.
Define an orthogonal matrix $Q$ as follows:
  \begin{align}\label{eq:2}
    Q=\left(
           \begin{matrix}
             Q_1   &  &  &  &  \\
                & Q_2 &  &  &  \\
                &  & \ddots &  &  \\
                &  &  & Q_k &  \\
                &  &  &  & I_{c} \\
           \end{matrix}
         \right).
  \end{align}
  According to the proof of \cite{GM}*{Theorem 2.2}, we know $Q^{-1}A(X)Q=A(X^{\pi})$.

  On the other hand, based on the partition $\pi$, we can write the
  degree matrices of $X$ and $X^{\pi}$ as follows:
  \begin{align*}
    D(X)=D(X^{\pi})=\left(
           \begin{matrix}
             g_1I_{c_1}   &  &  &  &  \\
                & g_2I_{c_2} &  &  &  \\
                &  & \ddots &  &  \\
                &  &  & g_kI_{c_k} &  \\
                &  &  &  & D(X[C])+(\sum_{i=1}^k \frac{c_i}{2})I_{c} \\
           \end{matrix}
         \right),
  \end{align*}
  where $g_i=\sum_{j=1}^kd_{ij}+\frac{c}{2}$  for $i\in \{1,2,\ldots,k\}$.
It is routine to check that $Q^{-1}D(X)Q=D(X^{\pi})$. Thus, $X$ and $X^{\pi}$
are degree similar.\qed
\begin{eg}
In Figure \ref{fig1}, $X_{1,2}$ can be obtained from $X_{1,1}$ by using local
switching. Let
  \begin{equation*}
    R_1=\left(
        \begin{matrix}
          \frac{1}{2}J_4-I_4 &  \\
           & I_6 \\
        \end{matrix}
      \right).
  \end{equation*}
Then $R_1^{-1}A(X_{1,1})R_1=A(X_{1,2})$ and $R_1^{-1}D(X_{1,1})R_1=D(X_{1,2})$.
Hence $X_{1,1}$ and $X_{1,2}$ are degree similar.
\begin{figure}[ht!]
  \centering
\includegraphics[width=100mm]{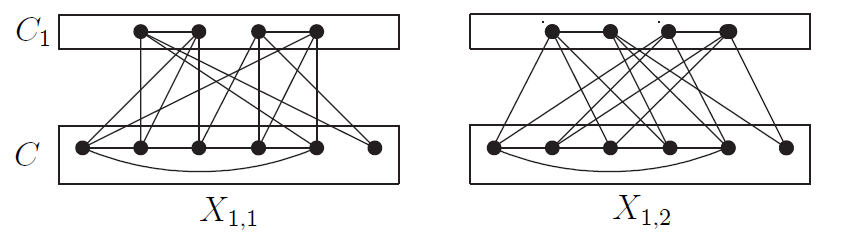}\\
  \caption{Local switching.}\label{fig1}
\end{figure}\qed
\end{eg}

\section{Joins and Products}

In this section, we will investigate how other standard graph theoretic operations
can be used to get examples of degree-similar graphs.

Let $X$ be a graph with two induced subgraphs $Y$ and $Z$ such that $V(X)$
is the disjoint union of $V(Y)$ and $V(Z)$, $E(X)$ is the disjoint union of $E(Y)$ and $E(Z)$.
We will say that $X$ is the \textsl{union} of $Y$ and $Z$, and denote it
by $Y\cup Z$.
The \textsl{join} of $Y$ and $Z$, written as $X\vee Y$, is the graph
obtained from $Y\cup Z$ by joining each vertex in $Y$ to each vertex in $Z$.

Now, we construct degree-similar graphs by using union and join operations.
The following results can be proved directly by the definition of degree-similar graphs.

\begin{lemma}
  Let $X$ and $Y$ be two connected degree-similar graphs. For any graph $H$, the following hold.
  \begin{enumerate}[(i)]
    \item $X\cup H$ and $Y\cup H$ are degree similar;
    \item If $X$ is regular, then $X\vee H$ and $Y\vee H$ are degree similar.
  \end{enumerate}
\end{lemma}

\begin{lemma}
  Let $X$ and $X^{\pi}$ be graphs defined in Lemma \ref{lem9} with vertex
  partition $\pi:=C_1\cup C_2\cup \cdots \cup C_k\cup C$, and let $Y$ be any graph.
  For convenience, put $C_{k+1}:=C$. In $X$ and $X^{\pi}$, if we add all
  edges between $Y$ and $C_{i_1}\cup C_{i_2}\cup \cdots\cup C_{i_l}$,
  where $1\leqslant i_1<i_2<\cdots<i_l\leqslant k+1$, then the two graphs obtained
  are degree similar.
\end{lemma}

\proof
 Denote by $\Gamma_1$ and $\Gamma_2$ the graphs obtained from $X$
 and $X^{\pi}$ by adding all edges between $Y$ and $C_{i_1}\cup C_{i_2}\cup \cdots\cup C_{i_l}$
 respectively. Assume $|V(X)|=n$ and $|V(Y)|=m$.
 Let $Q$ be the matrix defined in \eqref{eq:2}. Based on the proof of Lemma~\ref{lem9}, we know
  \[
    Q^{-1}A(X)Q=A(X^{\pi}),\quad Q^{-1}D(X)Q=D(X^{\pi}).
  \]
  Notice that
  \begin{equation*}
    A(\Gamma_1)=\left(
                   \begin{matrix}
                     A(X) & B \\
                     B^T & A(Y) \\
                   \end{matrix}
                 \right),\quad
    D(\Gamma_1)=\left(
                   \begin{matrix}
                     D(X)+mC &  \\
                      & D(Y)+(\sum_{j=1}^l|C_{i_j}|)I_{m} \\
                   \end{matrix}
                 \right).
  \end{equation*}
where $B=(b_{uv})_{n\times m}$ with $b_{uv}=1$ if $u\in C_{i_1}\cup C_{i_2}\cup \cdots\cup C_{i_l}$
and $b_{uv}=0$ otherwise; $C=(c_{uv})_{n\times n}$ is a diagonal matrix
with $c_{uu}=1$ if $u\in C_{i_1}\cup C_{i_2}\cup \cdots\cup C_{i_l}$ and $c_{uu}=0$ otherwise.

  Now, we define an orthogonal matrix as follows:
   \begin{align*}
    R=\left(
        \begin{matrix}
          Q &  \\
           & I_{m} \\
        \end{matrix}
      \right).
  \end{align*}
  Clearly, $R^T=R^{-1}=R$. By a direct calculation, one has
  $(\frac{2}{c_i}J_{c_i}-I_{c_i}){\bf 1}_{c_i}={\bf 1}_{c_i}$ for all $i\in \{1,\ldots,k\}$, where ${\bf 1}_{c_i}$ denotes the all 1s column vector of order $c_i$. Therefore, $QB=B$.  Hence
  \begin{align*}
    R^{-1}A(\Gamma_1)R&=\left(
                   \begin{matrix}
                     A(X^{\pi}) & QB \\
                     B^TQ & A(Y)\\
                   \end{matrix}
                 \right)=\left(
                   \begin{matrix}
                     A(X^{\pi}) & B \\
                   B^T & A(Y) \\
                   \end{matrix}
                 \right)=A(\Gamma_2),\\
    R^{-1}D(\Gamma_1)R&=\left(
                   \begin{matrix}
                     D(X^{\pi})+mC &  \\
                      & D(Y)+(\sum_{j=1}^l|C_{i_j}|)I_{m}\\
                   \end{matrix}
                 \right)=D(\Gamma_2).
  \end{align*}
  That is, $\Gamma_1$ and $\Gamma_2$ are degree similar.\qed

\begin{eg}
In Figure \ref{fig2}, $X_{2,i}$ is a graph obtained from $X_{1,i}$ by adding a
path $P_3$ and join it to all vertices of $C_1$ in $X_{1,i}$ for $i\in \{1,2\}$.

  \begin{figure}[ht!]
  \centering
  \includegraphics[width=100mm]{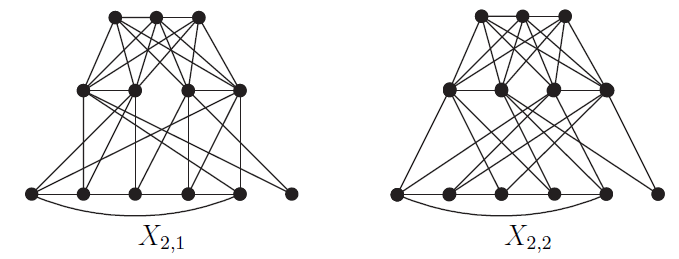}\\
  \caption{Joins of $P_3$ with a vertex subset $C_1$ of $X_{1,i}$.}\label{fig2}
\end{figure}

Define
  \begin{equation*}
    R_2=\left(
        \begin{matrix}
          \frac{1}{2}J_4-I_4 & & \\
           & I_6 & \\
           & & I_3\\
        \end{matrix}
      \right).
  \end{equation*}
It is routine to check that $R_2^{-1}A(X_{2,1})R_2=A(X_{2,2})$
and $R_2^{-1}D(X_{2,1})R_2=D(X_{2,2})$. Therefore, $X_{2,1}$
and $X_{2,2}$ are degree similar.\qed
\end{eg}

Recall that $\psi(X,t,\mu)=\psi(Y,t,\mu)$ is a necessary condition for two
graphs $X$ and $Y$ being degree similar. In fact, there is more we can say
about $\psi(X,t,\mu)$.
The following theorem is a generalization of Johnson and Newman's result
in \cite{MR565513}. For the detailed proof, one may see \cite{xiaohong}.

\begin{theorem}\label{thm3.1}
  Assume $X$ and $Y$ are graphs with $\psi(X,t,\mu)=\psi(Y,t,\mu)$ and $\psi(\overline{X},t,\mu)=\psi(\overline{Y},t,\mu)$, then there is an orthogonal matrix $Q$ such that
  \[
    Q^TA_\mu(X)Q=A_\mu(Y)\quad \text{and}\quad Q{\bf 1}={\bf 1}.
  \]
\end{theorem}

Based on the above theorem, by using join operator, we may obtain a family
of cospectral graphs with respect to the $A-\mu D$ matrix.

\begin{lemma}\label{lemma01}
If $X$ and $Y$ are graphs with
$\psi(X,t,\mu)=\psi(Y,t,\mu)$ and $\psi(\overline{X},t,\mu)=\psi(\overline{Y},t,\mu)$,
then $\psi(X\vee H,t,\mu)=\psi(Y\vee H,t,\mu)$
and $\psi(\overline{X\vee H},t,\mu)=\psi(\overline{Y\vee H},t,\mu)$,
where $H$ is any graph.
\end{lemma}

\proof
  Assume $X$ and $Y$ have $n$ vertices, $H$ has $m$ vertices.
  In view of Theorem \ref{thm3.1}, there is an orthogonal matrix $Q$ such that
  \[
    Q^TA_\mu(X)Q=A_\mu(Y)\quad \text{and}\quad Q{\bf 1}_n={\bf 1}_n.
  \]
  Then $Q^T{\bf 1}_n={\bf 1}_n$, which implies $Q^TJ_{n\times m}=J_{n\times m}$.
  Clearly, 
  \begin{equation*}
    A_\mu(X\vee H)=
     \left(
       \begin{matrix}
         A_\mu(X)-\mu nI_m & J_{n\times m} \\
         J_{m\times n} & A_\mu(H)-\mu mI_n \\
       \end{matrix}
     \right).
  \end{equation*}
  It is routine to check that
  \begin{equation*}
    \left(
       \begin{matrix}
         Q^T &  \\
          & I \\
       \end{matrix}
     \right)A_\mu(X\vee H)\left(
       \begin{matrix}
         Q &  \\
          & I \\
       \end{matrix}
     \right)=A_\mu(Y\vee H).
  \end{equation*}
  Hence $A_\mu(X\vee H)$ and $A_\mu(Y\vee H)$ are similar.
  It follows that $\psi(X\vee H,t,\mu)=\psi(Y\vee H,t,\mu)$.

Obviously, $\overline{X\vee H}=\overline{X}\cup \overline{H}$
and $\overline{Y\vee H}=\overline{Y}\cup \overline{H}$. Together with the
fact that $\psi(\overline{X},t,\mu)=\psi(\overline{Y},t,\mu)$,
then $\psi(\overline{X\vee H},t,\mu)=\psi(\overline{Y\vee H},t,\mu)$ holds obviously.\qed

Our next contribution involves constructions of degree-similar graphs
using alternative graph products, namely Cartesian product,  tensor product,
strong product and lexicographic product.

Let $A$ and $B$ be matrices of order $m\times n$ and $p\times q$,
respectively. The \textsl{Kr\"{o}necker product} of two matrices $A$ and $B$, denoted $A\otimes B$, is the $mp\times  nq$
block matrix $[a_{ij}B]$. It can be verified from the definition that
\[
(A\otimes B)(C\otimes D) = AC\otimes BD.
\]

Let $X$ and $Y$ be graphs with vertex sets $V(X)$ and $V(Y)$,
respectively. Put $n:=|V(X)|$ and $m:=|V(Y)|$.
The \textsl{Cartesian product} of $X$ and $Y$, denoted by $X\Box Y$,
is the graph defined as follows.
The vertex set of $X\Box Y$ is $V(X)\times V(Y)$.
The vertices $(u,v)$ and $(u', v')$ are adjacent
if either $u = u'$ and $v$ is adjacent to $v'$ in $Y$, or $v = v'$ and $u$
is adjacent to $u'$ in $X$. It is well known that
\[
  A(X\Box Y)=(A(X)\otimes I_m)+(I_n\otimes A(Y)),\ D(X\Box Y)=(D(X)\otimes I_m)+(I_n\otimes D(Y)).
\]

The \textsl{tensor product} of $X$ and $Y$, denoted by $X\otimes Y$,
is the graph defined as follows.
The vertex set of $X\otimes Y$ is $V(X)\times V(Y)$. The vertices $(u,v)$
and $(u', v')$ are adjacent
if $u$ is adjacent to $u'$ in $X$ and $v$ is adjacent to $v'$ in $Y$.
Notice that
\[
  A(X\otimes Y)=A(X)\otimes A(Y),\quad D(X\otimes Y)=D(X)\otimes D(Y).
\]

The \textsl{strong product} of $X$ and $Y$, denoted by $X\boxtimes Y$,
is the graph defined as follows.
The vertex set of $X\boxtimes Y$ is $V(X)\times V(Y)$. The vertices $(u,v)$
and $(u', v')$ are adjacent
if either $u = u'$ and $v$ is adjacent to $v'$ in $Y$, or $v = v'$ and $u$
is adjacent to $u'$ in $X$, or $u$ is adjacent to $u'$ in $X$ and $v$ is adjacent to $v'$ in $Y$. Then
\begin{align*}
  A(X\boxtimes Y)&=(A(X)\otimes I_m)+(I_n\otimes A(Y))+A(X)\otimes A(Y),\\
  D(X\boxtimes Y)&=(D(X)\otimes I_m)+(I_n\otimes D(Y))+D(X)\otimes D(Y).
\end{align*}

The \textsl{lexicographic product} of $X$ and $Y$, denoted by $X\odot Y$,
is the graph defined as follows.
The vertex set of $X\odot Y$ is $V(X)\times V(Y)$. The vertices $(u,v)$
and $(u', v')$ are adjacent
if either $u$ is adjacent to $u'$ in $X$, or $u = u'$ and $v$ is adjacent
to $v'$ in $Y$. It is easy to obtain
\begin{align*}
  A(X\odot Y)&=(A(X)\otimes J_m)+(I_n\otimes A(Y)),\\
  D(X\odot Y)&=(mD(X)\otimes I_m)+(I_n\otimes D(Y)).
\end{align*}

\begin{lemma}\label{lem5}
  If $X_1$ and $X_2$ are two degree-similar graphs with order $n$, $Y_1$
  and $Y_2$ are two degree-similar graphs with order $m$,
  then $X_1\Box Y_1$ and $X_2\Box Y_2$ (resp. $X_1\otimes Y_1$
  and $X_2\otimes Y_2$, $X_1\boxtimes Y_1$ and $X_2\boxtimes Y_2$) are degree similar.
  Furthermore, if $Y_1$ is connected, then $X_1\odot Y_1$ and $X_2\odot Y_2$
  are also degree similar.
\end{lemma}

\proof
  Here, we only prove $X_1\Box Y_1$ and $X_2\Box Y_2$ are degree similar,
  the remaining cases can be proved similarly.  By assumption, there exist two invertible real matrices $M_1$ and $M_2$ such that
\[
	M_1^{-1}A(X_1)M_1= A(X_2),\quad M_1^{-1}D(X_1)M_1 = D(X_2),
\]
and
\[
	M_2^{-1}A(Y_1)M_2= A(Y_2),\quad M_2^{-1}D(Y_1)M_2 = D(Y_2).
\]
Let $M=M_1\otimes M_2$. By applying the properties of Kr\"{o}necker product, one has
\begin{align*}
  M^{-1}A(X_1\Box Y_1)M&=(M_1\otimes M_2)^{-1}((A(X_1)\otimes I_m)+ (I_n\otimes A(Y_1)))(M_1\otimes M_2)\\
  &=(M_1^{-1}A(X_1)M_1\otimes M_2^{-1}I_mM_2)+(M_1^{-1}I_nM_1\otimes M_2^{-1}A(Y_1)M_2)\\
  &=(A(X_2)\otimes I_m)+(I_n\otimes A(Y_2))\\
  &=A(X_2\Box Y_2),
\end{align*}
and
\begin{align*}
  M^{-1}D(X_1\Box Y_1)M&=(M_1\otimes M_2)^{-1}((D(X_1)\otimes I_m)+ (I_n\otimes D(Y_1)))(M_1\otimes M_2)\\
  &=(M_1^{-1}D(X_1)M_1\otimes M_2^{-1}I_mM_2)+(M_1^{-1}I_nM_1\otimes M_2^{-1}D(Y_1)M_2)\\
  &=(D(X_2)\otimes I_m)+(I_n\otimes D(Y_2))\\
  &=D(X_2\Box Y_2).
\end{align*}
It follows that $X_1\Box Y_1$ and $X_2\Box Y_2$ are degree similar.\qed

\section{$k$-sum and rooted product}

In this section, we show how to use $k$-sum and rooted product of graphs to build
families of degree-similar graphs. The $k$-\textsl{sum} of graphs $X$ and $Y$ is
obtained by merging $k$ distinct vertices in $X$ with $k$ distinct vertices in $Y$.

\begin{lemma}\label{lem8}
  Let $X_1$ and $X_2$ be two degree-similar graphs, and $Y$ be an $n$-vertex
  graph with $\{w_1,\ldots,w_k\}\subseteq V(Y)$. Choose $u_1,\ldots,u_k\in V(X_1)$
  and $v_1,\ldots,v_k\in V(X_2)$ such that for $i\in \{1,\ldots,k\}$,
   \begin{enumerate}[(i)]
     \item the degree of $u_i$ (resp. $v_i$) is different with that of all other vertices in $X_1$ (resp. $X_2$);
     \item $d_{X_1}(u_i)=d_{X_2}(v_i)$;
     \item $X_1[\{u_1,\ldots,u_k\}]$ is connected.
   \end{enumerate}
   Denote by $\Gamma_1$ (resp. $\Gamma_2$) the $k$-sum of $X_1$
   (resp. $X_2$) and $Y$, which is obtained by merging $\{u_1,\ldots,u_k\}$
   of $X_1$ (resp. $\{v_1,\ldots,v_k\}$ of $X_2$) with $\{w_1,\ldots,w_k\}$
   of $Y$ in order. Then $\Gamma_1$ and $\Gamma_2$ are degree similar.
\end{lemma}

\proof
  Since $X_1$ and $X_2$ are two degree-similar graphs, there exists an
  invertible matrix $M$ such that
  \[
    M^{-1}A(X_1)M=A(X_2),\quad M^{-1}D(X_1)M=D(X_2).
  \]
  In view of Lemma \ref{lem1}, $M$ is block diagonal with respect to the
  partition $(V(X)\setminus \{u_1\})\cup \{u_1\}\cup \cdots \cup \{u_k\}$.
  Therefore, $M$ can be written as
  \begin{equation*}
    M=\left(
       \begin{matrix}
          M_1 & & & \\
          & a_1& & \\
          & & \ddots & \\
          & & & a_k\\
       \end{matrix}
     \right)
  \end{equation*}
for some invertible matrix $M_1$ and nonzero real numbers $a_1,\ldots,a_k$.
The adjacency matrix of $\Gamma_1$ is
\begin{equation*}
  A(\Gamma_1)=\left(
                \begin{matrix}
                  A(X_1\backslash \{u_1,\ldots,u_k\}) & {\bf x_1} & {\bf x_2} &\cdots & {\bf x_k} & {\bf 0} \\
                  {\bf x_1}^T & 0 & a_{12} &\cdots & a_{1k} & {\bf y_1}^T \\
                  {\bf x_2}^T & a_{21} & 0 &\cdots & a_{2k} & {\bf y_2}^T \\
                  \vdots & \vdots & \vdots &\ddots & \vdots & \vdots\\
                  {\bf x_k}^T & a_{k1} & a_{k2} &\cdots & 0 & {\bf y_k}^T \\
                  {\bf 0} & {\bf y_1} & {\bf y_2} &\cdots & {\bf y_k} & A(Y\backslash \{w_1,\ldots,w_k\}) \\
                \end{matrix}
              \right),
\end{equation*}
for some column vectors ${\bf x_i},\ {\bf y_i}$
and $a_{ij}\in \{0,1\}$ for $1\leqslant i,j\leqslant k$.
The degree matrix of $\Gamma_1$ is
\begin{equation*}
  D(\Gamma_1)=\left(
                \begin{matrix}
                  D(X_1)(\{u_1,\ldots,u_k\}) &  &  &  &  &  \\
                   & b_1 &  &  &  &  \\
                   &  & b_2 &  &  &  \\
                   &  &  & \ddots &  &  \\
                   &  &  &  & b_k &  \\
                   &  &  &  &  & D(Y)(\{w_1,\ldots,w_k\}) \\
                \end{matrix}
              \right),
\end{equation*}
where $b_i=d_{X_1}(u_i)+d_{Y}(w_i)-|\{j:u_iu_j\in E(X_1)\ \text{and}\ w_iw_j\in E(Y)\}|$.
Furthermore, the adjacency matrix and degree matrix of $\Gamma_2$ have the similar forms.

Recall that $M^{-1}A(X_1)M=A(X_2)$. Then $a_i^{-1}a_{ij}a_j=0$ if $a_{ij}=0$,
and $a_i^{-1}a_{ij}a_j=1$ if $a_{ij}=1$.
Therefore, $X_1[\{u_1,\ldots,u_k\}]\cong X_2[\{v_1,\ldots,v_k\}]$.
Together with $X_1[\{u_1,\ldots,u_k\}]$ is connected, one has $a_1=\cdots=a_k$. Let
\begin{equation*}
  Q=\left(
      \begin{matrix}
        M_1 &  &  \\
         & a_1I_k &  \\
         &  & a_1I_{n-k} \\
      \end{matrix}
    \right).
\end{equation*}
It is straightforward to check that
\[
  Q^{-1}A(\Gamma_1)Q=A(\Gamma_2),\quad Q^{-1}D(\Gamma_1)Q=D(\Gamma_2),
\]
i.e., $\Gamma_1$ and $\Gamma_2$ are degree similar.\qed

Taking a base graph $X$ and a sequence $\mathcal{Y}$ of rooted
graphs $Y_1,\ldots,Y_k$,
and then merge the $k$ roots of graphs in $\mathcal{Y}$ with $k$ distinct vertices $u_1,\ldots,u_k$ of $X$.
We refer to it as the \textsl{rooted product} of $X$ with $\mathcal{Y}$ at $u_1,\ldots,u_k$.
By a similar discussion as Lemma \ref{lem8}, we can get the following result.

\begin{lemma}
  Let $X_1$ and $X_2$ be two degree-similar graphs, and let $\mathcal{Y}=(Y_1,\ldots,Y_k)$
  be a sequence of rooted graphs. Choose $u_1,\ldots,u_k\in V(X_1)$
  and $v_1,\ldots,v_k\in V(X_2)$ such that for $i\in \{1,\ldots,k\}$,
   \begin{enumerate}[(i)]
     \item the degree of $u_i$ (resp. $v_i$) is different with that of all other vertices in $X_1$ (resp. $X_2$);
     \item $d_{X_1}(u_i)=d_{X_2}(v_i)$.
   \end{enumerate}
   Then the rooted product of $X_1$ with $\mathcal{Y}$ at $u_1,\ldots,u_k$ and the rooted product of $X_2$ with $\mathcal{Y}$ at $v_1,\ldots,v_k$ are degree similar.
\end{lemma}

\begin{eg}
In Figure \ref{fig3}, for $i\in \{1,2\}$, $X_{3,i}$ is the $2$-sum of $X_{1,i}+uv$
and a cycle $C_3$. Notice that the degree of $u$ (resp. $v$) is different with that
of all other vertices in $X_{1,i}+uv$. Let
  \begin{equation*}
    R_3=\left(
        \begin{matrix}
          \frac{1}{2}J_4-I_4 &  &\\
           & I_6 &\\
           & & 1\\
        \end{matrix}
      \right).
  \end{equation*}
One may check that $R_3^{-1}A(X_{3,1})R_3=A(X_{3,2})$ and $R_3^{-1}D(X_{3,1})R_3=D(X_{3,2})$. Hence $X_{3,1}$ and $X_{3,2}$ are degree similar.
\begin{figure}[ht!]
  \centering
  \includegraphics[width=100mm]{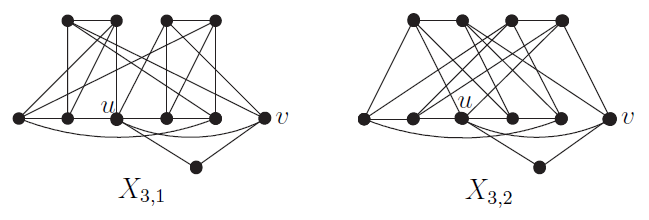}\\
  \caption{2-sum.}\label{fig3}
\end{figure}

In Figure~\ref{fig4}, for $i\in \{1,2\}$, $X_{4,i}$ is the rooted product of $X_{1,i}+uv$ with $(P_3,C_3)$ at $u,v$. Let
\begin{equation*}
  R_4=\left(
        \begin{matrix}
          \frac{1}{2}J_4-I_4 &  &  &  \\
           & I_6 &  &  \\
           &  & I_2 &  \\
           &  &  & I_2 \\
        \end{matrix}
      \right).
\end{equation*}
By a direct calculation, one has $R_4^{-1}A(X_{4,1})R_4=A(X_{4,2})$ and $R_4^{-1}D(X_{4,1})R_4=D(X_{4,2})$. Thus, $X_{4,1}$ and $X_{4,2}$ are degree similar.

  \begin{figure}[ht!]
  \centering
  \includegraphics[width=95mm]{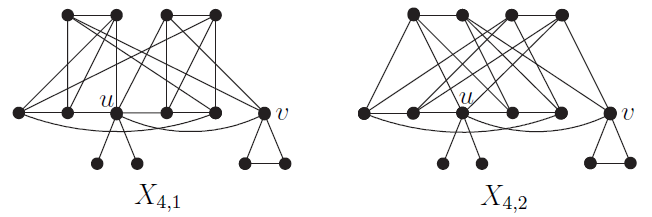}\\
  \caption{Rooted product.}\label{fig4}
\end{figure}\qed
\end{eg}

\section{Adding  or deleting  vertices}
Butler et al. \cite{butler2022complements} showed that: if $X_1$ and $X_2$
are two graphs with $\psi(X_1,t,\mu)=\psi(X_2,t,\mu)$, then the graphs
resulting from attaching an arbitrary rooted graph $Y$ to each vertex of $X_1$
and each vertex of $X_2$ will be cospectral with respect to the matrix $A-\mu D$.
Motivated by this, in this section, we construct degree-similar graphs by
adding or deleting vertices.

\begin{lemma}
  Let $X_1$ and $X_2$ be two degree-similar graphs. Assume $d_1, \ldots, d_k$
  are all different vertex degrees in $X_1$. For each $i\in \{1,\ldots,k\}$,
  attach $s_i$ pendant vertices to each vertex with degree $d_i$ in $X_1$
  and $X_2$ respectively, then the two graphs obtained  are degree similar.
\end{lemma}

\proof
  For convenience, denote by $\Gamma_1$ and $\Gamma_2$ the two graphs obtained from $X_1$ and $X_2$ respectively. Firstly, we reorder the vertices of $X_1$ and $X_2$ such that their degree matrices can be written as
  \begin{equation*}
    D(X_1)=D(X_2)=\left(
                    \begin{matrix}
                      d_1I_{n_1} &  &  &  \\
                       & d_2I_{n_2} &  &  \\
                       &  & \ddots &  \\
                       &  &  & d_kI_{n_k} \\
                    \end{matrix}
                  \right),
  \end{equation*}
  where $n_i$ is the number of vertices in $V_i:=\{v: d_{X_1}(v)=d_i\}$
  for $i\in \{1,\ldots,k\}$.

  Since $X_1$ and $X_2$ are degree similar, there exists an
  invertible real matrix $M$ such that
  \[
    M^{-1}A(X_1)M=A(X_2),\quad M^{-1}D(X_1)M=D(X_2).
  \]
  In view of Lemma \ref{lem1}, $M$ is block diagonal with respect to the
  partition $V_1\cup \cdots \cup V_k$. That is,
  \begin{equation*}
    M=\left(
                    \begin{matrix}
                      M_1 &  &  &  \\
                       & M_2 &  &  \\
                       &  & \ddots &  \\
                       &  &  & M_k \\
                    \end{matrix}
                  \right)
  \end{equation*}
  for some invertible matrices $M_1,\ldots,M_k$. For $i\in \{1,\ldots,k\}$,
  assume
  \[
  V_i=\{v_{n_1+\cdots+n_{i-1}+1},v_{n_1+\cdots+n_{i-1}+2},\ldots,v_{n_1+\cdots+n_{i-1}+n_i}\},
  \]
  and for $j\in \{1,\ldots,n_i\}$, assume $\{u^i_{(j-1)s_i+1},\ldots,u^i_{js_i}\}$ are all pendant vertices attached to
  $v_{n_1+\cdots+n_{i-1}+j}$. Next, partition the vertex set of $\Gamma_1$ as follows:
  \begin{align*}
  &\{u^1_1,u^1_{s_1+1},\ldots,u^1_{(n_1-1)s_1+1}\}
  \cup \{u^1_2,u^1_{s_1+2},\ldots,u^1_{(n_1-1)s_1+2}\}
  \cup \cdots\cup \{u^1_{s_1},u^1_{2s_1},\ldots,u^1_{n_1s_1}\}\\
  \cup& \{u^2_1,u^2_{s_2+1},\ldots,u^2_{(n_2-1)s_2+1}\}
  \cup \{u^2_2,u^2_{s_2+2},\ldots,u^2_{(n_2-1)s_2+2}\}
  \cup \cdots\cup \{u^2_{s_2},u^2_{2s_2},\ldots,u^2_{n_2s_2}\}\\
  \cup &\cdots\\
  \cup& \{u^k_1,u^k_{s_k+1},\ldots,u^k_{(n_k-1)s_k+1}\}
  \cup \{u^k_2,u^k_{s_k+2},\ldots,u^k_{(n_k-1)s_k+2}\}
  \cup \cdots\cup \{u^k_{s_k},u^k_{2s_k},\ldots,u^k_{n_ks_k}\}\\
  \cup & V_1\cup V_2\cup \cdots\cup V_k.
  \end{align*}

  It is routine to check that for $i\in \{1,2\}$, the adjacency matrix of $\Gamma_i$ is
  \begin{equation*}
    A(\Gamma_i)= \begin{pmatrix}
    {\bf 0} &
      \begin{matrix}
        {\bf 1}_{s_1}\otimes I_{n_1} & {\bf 0} & \cdots & {\bf 0}\\
       {\bf 0} & {\bf 1}_{s_2}\otimes I_{n_2}& \cdots & {\bf 0}\\
       \vdots &\vdots & \vdots &\vdots \\
       {\bf 0}&{\bf 0}&\cdots& {\bf 1}_{s_k}\otimes I_{n_k}\\
      \end{matrix}\\
      \begin{matrix}
        {\bf 1}_{s_1}^T\otimes I_{n_1} & {\bf 0} & \cdots & {\bf 0}\\
       {\bf 0} & {\bf 1}_{s_2}^T\otimes I_{n_2}& \cdots & {\bf 0}\\
       \vdots &\vdots & \vdots &\vdots \\
       {\bf 0}&{\bf 0}&\cdots &{\bf 1}_{s_k}^T\otimes I_{n_k}\\
      \end{matrix} & A(X_i)
    \end{pmatrix},
  \end{equation*}
  the degree matrix of $\Gamma_i$ is
  \begin{equation*}
    D(\Gamma_i)=\left(
                  \begin{array}{ccccc}
                    I_{n_1s_1+\cdots+n_ks_k} &  &  &  &  \\
                     & (d_1+s_1)I_{n_1} &  &  &  \\
                     &  & (d_2+s_2)I_{n_2} &  &  \\
                     &  &  & \ddots &  \\
                     &  &  &  & (d_k+s_k)I_{n_k} \\
                  \end{array}
                \right).
  \end{equation*}
  Let $Q$ be a matrix defined by
  \begin{equation*}
    Q=\left(
                  \begin{array}{ccccc}
                    I_{s_1}\otimes M_1 &  &  &  &  \\
                     & I_{s_2}\otimes M_2 &  &  &  \\
                     &  & \ddots &  &  \\
                     &  &  & I_{s_k}\otimes M_k &  \\
                     &  &  &  & M \\
                  \end{array}
                \right).
  \end{equation*}
  By a direct calculation, one has
   \[
  Q^{-1}A(\Gamma_1)Q=A(\Gamma_2),\quad Q^{-1}D(\Gamma_1)Q=D(\Gamma_2),
\]
i.e., $\Gamma_1$ and $\Gamma_2$ are degree similar.\qed

\begin{lemma}
  Let $X_1$ and $X_2$ be two degree-similar graphs. Assume $d_1, \ldots, d_k$
  are all different vertex degrees in $X_1$, and $V_i=\{v: d_{X_1}(v)=d_i\}$ with
  cardinality $n_i$. For each $i\in \{1,\ldots,l\}$ with $l\leqslant k$, add $n_i$ isolated
  vertices and join each of them to all vertices with degree $d_i$ in $X_1$ and $X_2$
  respectively, then the obtained graphs are degree similar.
\end{lemma}

\proof
  Denote by $\Gamma_1$ and $\Gamma_2$ the two graphs obtained from $X_1$ and $X_2$.
  Firstly, we reorder the vertices of $X_1$ and $X_2$ such that their degree matrices can be written as
  \begin{equation*}
    D(X_1)=D(X_2)=\left(
                    \begin{matrix}
                      d_1I_{n_1} &  &  &  \\
                       & d_2I_{n_2} &  &  \\
                       &  & \ddots &  \\
                       &  &  & d_kI_{n_k} \\
                    \end{matrix}
                  \right).
  \end{equation*}

  Since $X_1$ and $X_2$ are degree similar, there exists an
  invertible real matrix $M$ such that
  \[
    M^{-1}A(X_1)M=A(X_2),\quad M^{-1}D(X_1)M=D(X_2).
  \]
  In view of Lemma \ref{lem1}, $M$ is block diagonal with respect to the
  partition $V_1\cup \cdots \cup V_k$. That is,
  \begin{equation*}
    M=\left(
                    \begin{matrix}
                      M_1 &  &  &  \\
                       & M_2 &  &  \\
                       &  & \ddots &  \\
                       &  &  & M_k \\
                    \end{matrix}
                  \right)
  \end{equation*}
  for some invertible matrices $M_1,\ldots,M_k$.
  For $i\in \{1,\ldots,l\}$, assume
  \[
  	V_i = \{v_{n_1+\cdots+n_{i-1}+1},v_{n_1+\cdots+n_{i-1}+2},
  		\ldots,v_{n_1+\cdots+n_{i-1}+n_i}\},
  \]
  and assume $\{u^i_{1},\ldots,u^i_{n_i}\}$ are all added vertices that are adjacent to vertices
  in $V_i$.
  Partition the vertex set of $\Gamma_1$ as follows:
  \begin{align*}
  \{u^1_{1},\ldots,u^1_{n_1}\}
  \cup \{u^2_{1},\ldots,u^2_{n_2}\}
  \cup \cdots\cup \{u^l_{1},\ldots,u^l_{n_l}\}  \cup V_1\cup V_2\cup \cdots\cup V_l\cup (V_{l+1}\cup \cdots\cup V_k).
  \end{align*}
 For $i\in \{1,2\}$, it is routine to check that the adjacency matrix
 of $\Gamma_i$ is
  \begin{equation*}
    A(\Gamma_i)=
                      \begin{pmatrix}
                    {\bf 0} &  \begin{matrix}
                             J_{n_1} &  &  & {\bf 0} \\
                              & \ddots &   &\vdots \\
                              &  & J_{n_l} & {\bf 0} \\
                           \end{matrix}
                     \\
                           \begin{matrix}
                             J_{n_1} &  &    \\
                              & \ddots &    \\
                              &  & J_{n_l}   \\
                             {\bf 0} & \cdots &   {\bf 0} \\
                           \end{matrix} & A(X_i)\\
                  \end{pmatrix},
  \end{equation*}
  and the degree matrix of $\Gamma_i$ is
  \begin{equation*}
    D(\Gamma_i)=\left(
                  \begin{matrix}
                    n_1I_{n_1} &  &  &  &  &  &  &  &  \\
                     & \ddots &  &  &  &  &  &  &  \\
                     &  & n_lI_{n_l} &  &  &  &  &  &  \\
                     &  &  & (d_1+n_1)I_{n_1} &  &  &  &  &  \\
                     &  &  &  & \ddots &  &  &  &  \\
                     &  &  &  &  & (d_l+n_l)I_{n_l} &  &  &  \\
                     &  &  &  &  &  & d_{l+1}I_{n_{l+1}} &  &  \\
                     &  &  &  &  &  &  & \ddots &  \\
                     &  &  &  &  &  &  &  & d_kI_{n_k} \\
                  \end{matrix}
                \right).
  \end{equation*}
  Let
  \begin{equation*}
    Q=\left(
        \begin{matrix}
          M_1 &  &  &  \\
           & \ddots &  &  \\
           &  & M_l &  \\
           &  &  & M \\
        \end{matrix}
      \right).
  \end{equation*}
  Together with \eqref{eq:1}, and by a direct calculation, one has
   \[
  Q^{-1}A(\Gamma_1)Q=A(\Gamma_2),\quad Q^{-1}D(\Gamma_1)Q=D(\Gamma_2),
\]
i.e., $\Gamma_1$ and $\Gamma_2$ are degree similar.\qed

Finally, we construct degree-similar graphs by deleting a vertex.

\begin{lemma}
  Let $X_1$ and $X_2$ be two degree-similar graphs. For $i\in \{1,2\}$, choose $u_i\in V(X_i)$ such that
  \begin{enumerate}[(i)]
    \item the degree of $u_i$ in $X_i$ is different with that of all other vertices in $X_i$;
    \item $d_{X_1}(u_1)=d_{X_2}(u_2)$;
    \item $w\in N_{X_1}(u_1)$ implies $\{w':d_{X_1}(w')=d_{X_1}(w)\}\subseteq N_{X_1}(u_1)$.
  \end{enumerate}
  Then $X_1\backslash u_1$ and $X_2\backslash u_2$ are degree similar.
\end{lemma}

\proof
  Since $X_1$ and $X_2$ are two degree-similar graphs, there exists an invertible real  matrix $M$ such that
  \[
    M^{-1}A(X_1)M=A(X_2),\quad M^{-1}D(X_1)M=D(X_2).
  \]
  In view of Lemma \ref{lem1}, $M$ is block diagonal with respect to the partition $\{u_1\}\cup N_{X_1}(u_1) \cup (V(X_1)\setminus N_{X_1}[u_1])$. That is, $M$ can be written as follows:
  \begin{equation*}
    M=\left(
        \begin{matrix}
          a &  &  \\
           & M_1 &  \\
           &  & M_2 \\
        \end{matrix}
      \right),
  \end{equation*}
  for some invertible matrices $M_1,\,M_2$ and a nonzero real number $a$.   Furthermore, we can partition $A(X_1)$ and $D(X_1)$ as follows:
  \begin{equation*}
    A(X_1)=\left(
             \begin{matrix}
               0 & {\bf 1}^T & {\bf 0} \\
               {\bf 1} & A_{11} & A_{12} \\
               {\bf 0} & A_{21} & A_{22} \\
             \end{matrix}
           \right),\quad
    D(X_1)=\left(
             \begin{matrix}
               d_{X_1}(u_1) &  &  \\
                & D_1 &  \\
                &  & D_2 \\
             \end{matrix}
           \right).
  \end{equation*}
  Notice that the adjacency matrix and the degree matrix of $X_1\backslash u_1$ are
  \begin{equation*}
    A(X_1\backslash u_1)=\left(
                 \begin{matrix}
                   A_{11} & A_{12} \\
                   A_{21} & A_{22} \\
                 \end{matrix}
               \right),\quad
    D(X_1\backslash u_1)=\left(
                 \begin{matrix}
                   D_1-I &  \\
                    & D_2 \\
                 \end{matrix}
               \right).
  \end{equation*}
  Recall that $M^{-1}A(X_1)M=A(X_2)$ and $M^{-1}D(X_1)M=D(X_2)$, then
  \begin{equation*}
    A(X_2)=\left(
             \begin{matrix}
               0 & a^{-1}{\bf 1}^TM_1 & {\bf 0} \\
               M_1^{-1}{\bf 1}a &  M_1^{-1}A_{11} M_1 &  M_1^{-1}A_{12} M_2 \\
               {\bf 0} &  M_2^{-1}A_{21} M_1 &  M_2^{-1}A_{22} M_2 \\
             \end{matrix}
           \right),
  \end{equation*}
  and
  \begin{equation*}
           D(X_2)=\left(
             \begin{matrix}
               d_{X_1}(u_1) &  &  \\
                & M_1^{-1}D_1M_1 &  \\
                &  & M_2^{-1}D_2M_2 \\
             \end{matrix}
           \right).
  \end{equation*}
Clearly, the first row of $D(X_2)$ is indexed by $u_2$. Together
with Item (ii), we know each entry of $a^{-1}{\bf 1}^TM_1$ is equal to $1$.
Thus, $A(X_2)$ is partitioned according to
$\{u_2\}\cup N_{X_2}(u_2) \cup (V(X_2)\setminus N_{X_2}[u_2])$. Let $Q=\diag(M_1,M_2)$.
Then
\[
  Q^{-1}A(X_1\backslash u_1)Q=A(X_2\backslash u_2),\quad Q^{-1}D(X_1\backslash u_1)Q=D(X_2\backslash u_2),
\]
i.e., $X_1 \backslash u_1$ and $X_2\backslash u_2$ are degree similar.\qed

\begin{eg}

In Figure~\ref{fig6}, for $i\in \{1,2\}$, $X_{5,i}$ is a graph obtained
from $X_{1,i}+uv$ by adding one pendant vertex to each vertex with degree $4$,
adding two pendant vertices to $u$ and adding three pendant vertices to $v$. Let
\begin{equation*}
    R_5=\left(
        \begin{matrix}
          \frac{1}{2}J_4-I_4 & & & & &\\
          & I_4 & & & &\\
          & & I_2 & & &\\
          & &  &I_3 & &\\
          & &  & & \frac{1}{2}J_4-I_4&\\
          & & & & & I_6\\
        \end{matrix}
      \right).
  \end{equation*}
It is straightforward to check that $R_5^{-1}A(X_{5, 1})R_5=A(X_{5, 2})$
and $R_5^{-1}D(X_{5, 1})R_5=D(X_{5, 2})$. Hence $X_{5, 1}$ and $X_{5, 2}$
are degree similar.
\begin{figure}[ht!]
  \centering
  \includegraphics[width=95mm]{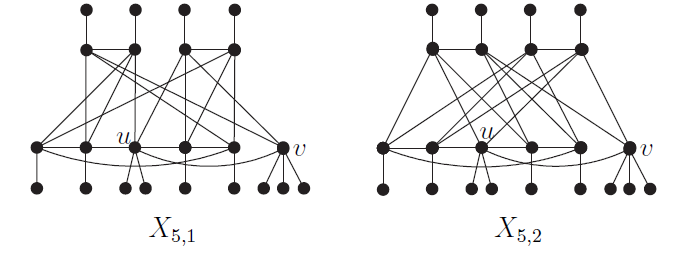}\\
  \caption{Adding pendant vertices.}\label{fig6}
\end{figure}

  In Figure~\ref{fig7}, by using Lemma \ref{lem9}, we know $X_{6,1}\backslash\{w_1,w_2\}$
  and $X_{6, 2}\backslash\{w_1,w_2\}$ are degree similar. Let
  \begin{equation*}
    R_6=\left(
        \begin{matrix}
          I_2 & & \\
          & \frac{1}{2}J_4-I_4 & \\
          & & I_6\\
        \end{matrix}
      \right).
  \end{equation*}
By a direct calculation, one has $R_6^{-1}A(X_{6, 1})R_6=A(X_{6, 2})$
and $R_6^{-1}D(X_{6, 1})R_6=D(X_{6, 2})$. Then $X_{6, 1}$ and $X_{6, 2}$
are degree similar.
\begin{figure}
  \centering
 \includegraphics[width=100mm]{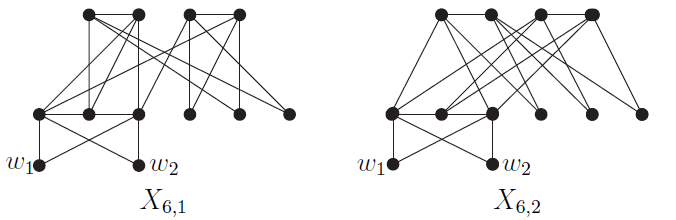}\\
  \caption{Adding complete graphs.}\label{fig7}
\end{figure}

In Figure \ref{fig5}, for $i\in \{1,2\}$, $X_{7,i}$ can be viewed as a graph
obtained from $X_{2,i}$ by deleting the unique vertex with degree $6$.
Notice that the neighborhood of this vertex contains all vertices with
degree $5$ and $7$ in $X_{2,i}$. Let
\begin{equation*}
    R_7=\left(
        \begin{matrix}
          \frac{1}{2}J_4-I_4 & & \\
           & I_6 & \\
           & & I_2\\
        \end{matrix}
      \right).
  \end{equation*}
Then $R_7^{-1}A(X_{7,1})R_7=A(X_{7,2})$ and $R_7^{-1}D(X_{7,1})R_7=D(X_{7,2})$.
Hence $X_{7,1}$ and $X_{7,2}$ are degree similar.
\begin{figure}[ht!]
  \centering
  \includegraphics[width=100mm]{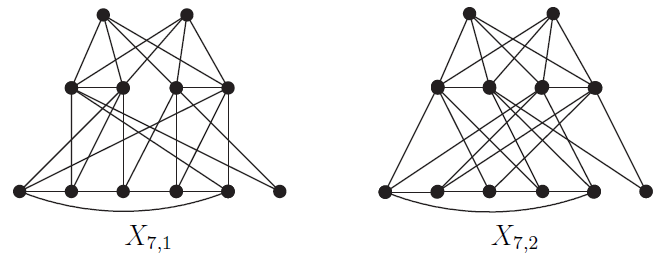}\\
  \caption{Deleting a vertex.}\label{fig5}
\end{figure}\qed
  \end{eg}

\section{Similarity}

Recall that
\[
	\psi(X,t,\mu)= \det(tI - (A-\mu D)).
\]
Here we view $A-\mu D$ as a matrix over the field of rational
functions $\rats(\mu)$, and then $\ptmu{X}$ is the characteristic polynomial
of $A-\mu D$.

We are accustomed to the fact that symmetric matrices
over $\re$ are similar if and only if their characteristic polynomials are
equal. In general though, equality of characteristic polynomials
is not enough to ensure similarity. Instead we have the following.
(For a proof and more details, see Lancaster and Tismenetsky
\cite{lancaster1985theory}*{Theorem~7.6.1} or Friedland
\cite{friedland2015matrices}*{Theorem~2.1.4}.)

\begin{theorem}\label{th9.1}
	Two matrices $B_1$ and $B_2$ over the field $\fld$ are similar
	if and only if the matrices $tI-B_1$ and $tI-B_2$ have the same
	Smith normal form.\qed
\end{theorem}

If $A$ and $D$ are integer matrices, then the Smith normal form of $tI-(A-\mu D)$
is determined by the determinants of submatrices of $tI-(A-\mu D)$, and
by the greatest common divisors of sets of these polynomials. Thus all
calculations are carried out in the principal ideal domain $\rats(\mu)[t]$.
This yields the following:

\begin{lemma}
	Matrices $A_1-\mu D_1$ and $A_2-\mu D_2$ are similar over $\re(\mu)$
	if and only if they are similar over $\rats(\mu)$.\qed
\end{lemma}

If $X_1$ and $X_2$ are degree-similar graphs with adjacency matrices
$A_1$, $A_2$ and degree matrices $D_1$, $D_2$,
then in view of Theorem~\ref{th9.1}, the
Smith normal forms of $tI-(A_1-\mu D_1)$ and $tI-(A_2-\mu D_2)$ are equal.
For the converse, if $tI-(A_1-\mu D_1)$ and $tI-(A_2-\mu D_2)$
have the same Smith normal form, then $A_1-\mu D_1$
and $A_2-\mu D_2$ are similar over $\rats(\mu)$. Furthermore, we have the
following result.

\begin{lemma}\label{lem9.2}
	Let $X_1$ and $X_2$ be graphs with adjacency matrices
	$A_1$, $A_2$ and degree matrices $D_1$, $D_2$ respectively.
	If $A_1-\mu D_1$ and $A_2-\mu D_2$ are similar over $\rats(\mu)$,
	then $A_1$ and $A_2$ are similar over $\rats$, as are $D_1$ and $D_2$.
\end{lemma}

\proof
Assume $A_1-\mu D_1$ and $A-\mu D_2$ are similar over $\rats(\mu)$.
Then there is a matrix $M(\mu)$, with entries rational functions in $\mu$,
such that
\begin{equation}
	\label{eq:MmuA1}
	M(\mu)^{-1}(A_1+\mu D_1)M(\mu) = A_2+\mu D_2.
\end{equation}
The set of poles of the entries of $M$ is finite, and therefore
for all sufficiently large rational numbers $\ga$, the real matrix $M(\ga)$
is invertible. Hence it follows from Equation~\eqref{eq:MmuA1}
that $D_1$ and $D_2$ are similar.

Next, there is a sequence of rational numbers $(\ga_i)_{i\ge0}$ converging
to zero such that $M(\ga_i)$ is defined and invertible, and it follows
that $A_1$ and $A_2$ are similar.\qed

\section{Problems}

In Sections 4-8, we provide a number of constructions of pairs of (non-isomorphic)
degree-similar graphs. It will be interesting to get more degree-similar graphs. In particular, the result in \cite{mckay1977spectral}*{Theorem 5.3} implies that two trees are degree similar if and only if they are isomorphic. A \textsl{unicyclic graph} can be viewed as a graph obtained from a tree by adding one edge. So, we present the first problem:
\begin{problem}
  Find more degree-similar graphs. In particular, are there non-isomorphic degree-similar unicyclic graphs?
\end{problem}

Based on Lemma \ref{lem9.2}, we know that if $tI-(A_1-\mu D_1)$
and $tI-(A_2-\mu D_2)$ have the same Smith normal form,
then $A_1$ and $A_2$ are similar over $\rats$, as are $D_1$ and $D_2$.
Then, a natural problem arises.

\begin{problem}
  Let $X$ and $Y$ be two graphs. Assume that $tI-(A(X)-\mu D(X))$
  and $tI-(A(Y)-\mu D(Y))$ have the same Smith normal form.
  Are $X$ and $Y$ degree similar?
\end{problem}

For a graph $X$ and an edge $e\in E(X)$, denote by $X\backslash e$
the graph obtained from $X$
by deleting the edge $e$. In \cite{xiaohong}, the authors showed that if $X$
is a strongly regular graph, then for any two edges $e$ and $f$ of $X$,
the graphs $X \backslash e$ and $X \backslash f$ are cospectral
with cospectral complements, with respect to the adjacency,
Laplacian, unsigned Laplacian and normalized Laplacian matrices.
Motivated by this, we consider a more general problem.

\begin{problem}
  Let $X$ be a strongly regular graph with two different edges $e$ and $f$.
  Are $X \backslash e$ and $X \backslash f$ degree similar?
\end{problem}


\bibliographystyle{abbrv}
\bibliography{degcosp}
\end{document}